\documentclass[12pt,draft]{article}

\usepackage{amsthm,amssymb,amsmath,amscd,amstext,mathrsfs}
\usepackage{latexsym, mathptmx}

\newtheorem{theorem}{Theorem}[section]
\newtheorem{definition}{Definition}[section]
\newtheorem{lemma}{Lemma}[section]

\newtheorem*{remark}{{\it Remark}}

\newcommand{\nc}{\newcommand}
\nc{\C}{{\mathbb C}}
\nc{\R}{{\mathbb R}}
\nc{\HH}{{\mathbb H}}
\nc{\Z}{{\mathbb Z}}
\nc{\N}{{\mathbb N}}   
\nc{\dd}{{\rm d}}
\nc{\DD}{{\rm D}}
\nc{\ii}{{\bf i}}
\nc{\ca}{{\mathscr A}}
\nc{\cb}{{\mathscr B}}
\nc{\cg}{{\mathscr G}}
\nc{\ch}{{\mathscr H}}
\nc{\cm}{{\mathscr M}}
\nc{\co}{{\mathscr O}}

\begin{document}

\title{On the energy spectrum of Yang--Mills instantons\\ 
over asymptotically locally flat spaces}

\author{G\'abor Etesi\\
\small{{\it Department of Geometry, Mathematical Institute, Faculty of
Science,}}\\
\small{{\it Budapest University of Technology and Economics,}}\\
\small{{\it Egry J. u. 1, H \'ep., H-1111 Budapest, Hungary
\footnote{{\tt etesi@math.bme.hu}}}}} 

\maketitle

\pagestyle{myheadings}
\markright{G.Etesi: On the energy spectrum of Yang--Mills instantons 
over ALF spaces}

\thispagestyle{empty}

\begin{abstract}
In this paper we prove that over an asymptotically locally flat 
(ALF) Riemannian four-manifold the energy of an ``admissible'' ${\rm SU}(2)$ 
Yang--Mills instanton is always integer. This result sharpens the previously 
known energy identity for such Yang--Mills instantons over ALF geometries. 
Furthermore we demonstrate that this statement continues to hold for 
the larger gauge group ${\rm U}(2)$.

Finally we make the observation that there might be a natural relationship  
between $4$ dimensional Yang--Mills theory over an ALF space and $2$ 
dimensional conformal field theory. This would provide a further 
support for the existence of a similar correspondence investigated by 
several authors recently.

\end{abstract}
\centerline{AMS Classification: Primary: 81T13; Secondary: 53C07, 58D27}
\centerline{Keywords: {\it ALF gravitational instanton; Yang--Mills 
instanton; energy}}


\section{Introduction}\label{one}


{\it Asymptotically locally flat} (ALF) spaces are non-compact 
complete Riemannian $4$-manifolds including mathematically as 
well as physically important spaces such as the Riemannian 
Schwarz\-schild and Kerr geometries \cite{egu-gil-han} and interesting 
hyper-K\"ahler examples like the flat $\R^3\times S^1$, the 
$A_k$ ALF (or multi-Taub--NUT or ALF Gibbons--Hawking) geometries 
\cite{gib-haw}, the Atiyah--Hitchin manifold 
\cite{ati-hit} and its cousins the so-called $D_k$ ALF spaces \cite{che-hit}. 

Recently there has been some effort to understand Yang--Mills instantons over 
these spaces from both mathematical (e.g. \cite{cha-hur1, cha-hur2, ete-hau1, 
ete-hau2, ete-hau3, ete-jar, ete-sza, nye}) and physical (e.g. 
\cite{bru-nog-van1, bru-nog-van2, che1, che2, che3, lee-yi, mur-voz, nak-sak, 
war, wit2}) sides. A central question is whether or not their moduli spaces are 
finite dimensional manifolds (with mild isolated singularities). Experienced 
with the compact case an expected condition is the discreteness of their energy 
spectrum. For instance in \cite{ete-jar} natural asymptotical analytical 
conditions on these solutions have been imposed in this spirit: the energy 
spectrum of these ``admissible'' Yang--Mills instantons is characterized by 
Chern--Simons invariants of the infinitely distant boundary of the ALF space. 
Hence the energy spectrum is discrete but in principle can contain fractional 
values correspondig to smooth Yang--Mills instantons with non-trivial holonomy 
at infinity. Over an $A_k$ ALE (or multi-Eguchi--Hanson or ALE 
Gibbons--Hawking) space smooth irreducible ${\rm SU}(2)$ Yang--Mills 
instanton solutions forming nice moduli spaces indeed can have fractional 
energy (cf. \cite{bia-fuc-ros-mar} or \cite[Theorem 4.1]{ete-hau3}). 
On the contrary over ALF spaces---although finite energy solutions are 
relatively easy to find---no ``admissible'' Yang--Mills instantons of 
fractional energy are known to exist (e.g. \cite{ete-hau1, ete-hau3, ete-sza}). 
Meanwhile other known finite energy solutions 
(e.g. \cite{bru-nog-van1, bru-nog-van2, cha-hur1, cha-hur2, che3, der, 
mos-tav}) with integer or fractional energy in principle may fit into larger 
continuous energy families, as an example in Sect. \ref{four} here shows.

The paper is organized as follows. In Sect. \ref{two} we argue why 
``admissibility'' (cf. Definition \ref{feltetelek} here or
\cite[Definition 2.1]{ete-sza}) is a quite weak and good condition to impose on 
${\rm SU}(2)$ Yang--Mills instantons in the ALF scenario: we shall see shortly 
that it is a sharp but easily satisfied condition moreover can be used in 
{\it all} known classes of ALF spaces (both the K\"ahler and non-K\"ahler 
examples) to rule out ${\rm SU}(2)$ Yang--Mills instantons with continuous 
energy. ``Admissibility'' rests on two important technicalities: 
the Hausel--Hunsicker--Mazzeo \cite{hau-hun-maz} compactification of 
an ALF space and the codimension $2$ singularity removal theorem of 
Sibner--Sibner \cite{sib-sib} and R\aa de \cite{rad}.

In Sect. \ref{three} we come to our main result and demonstrate that 
``admissibility'' excludes not only continuous but even fractional energy 
Yang--Mills instantons i.e., {\it over an ALF space the energy of a smooth 
``admissible'' ${\rm SU}(2)$ Yang--Mills instanton is always integer} although 
its associated holonomy at infinity is not necessarily trivial (cf. 
Theorem \ref{fotetel} here). 

A straightforward generalization in Sect. \ref{four} shows that our 
result remains valid for the larger gauge group ${\rm SU}(2)\times {\rm U}(1)$ 
and---under some topological conditions on the underlying vector 
bundle---also for ${\rm U}(2)$ although when deriving these results some care 
is needed due to the existence of non-topological ${\rm U}(1)$ Yang--Mills 
instantons (cf. Theorem \ref{U2} here). 

Finally in Sect. \ref{five} we speculate if the existence of the 
aforementioned compactification of an ALF space with its pleasant 
properties in Yang--Mills theory might indicate an intrinsic 
AGT-like relationship at the quantum level between $4d$ YM theory and $2d$ CFT 
as it has been suggested recently by several authors from a different 
angle, cf. e.g. \cite{ald-gai-tac, alf-tar, bel-ber-fei-lit-tar, tac}.


\section{An overview of Yang--Mills theory over ALF spaces}
\label{two}


So let $(M,g)$ be an ALF space as defined in \cite{ete-jar, hau-hun-maz}. 
Topologically, an ALF space (with a single end) admits a decomposition 
$M=K\cup W$ where $K$ is a compact interior space and $W$ is an end or neck 
homeomorphic to $N\times\R^+$ where 
\begin{equation}
\pi :N\longrightarrow B_{+\infty}
\label{fibralas}
\end{equation}
is a connected, compact, oriented three-manifold fibered over a connected 
compact surface $B_{+\infty}$ with circle fibers $F\cong S^1$. Consider 
$W\cong N\times\R^+$ and denote by $r\in [0, +\infty )$ the radial coordinate 
parameterizing $\R^+$. Regarding the complete Riemannian metric $g$ 
there exists a diffeomorphism $\phi :N\times\R^+\rightarrow W$ such that
\[\phi^*(g\vert_W)=\dd r^2+r^2(\pi^*g_{B_{+\infty}})'+h'_F\]
where $g_{B_{+\infty}}$ is a smooth metric on $B_{+\infty}$, $h_F$ is a
symmetric 2-tensor on $N$ which restricts to a metric along the
fibers $F\cong S^1$ and $(\pi^*g_{B_{+\infty}})'$ as well as $h'_F$ are 
smooth non-vanishing $O(1)$ extensions of $\pi^*g_{B_{+\infty}}$ and $h_F$ 
over $W$, respectively. Furthermore, we impose that the curvature $R_g$ of 
$g$ decays like
\[\vert\phi^*(R_g\vert_W)\vert\sim O(r^{-3}).\]
Here $R_g$ is regarded as a map $R_g:C^\infty (\Lambda^2M)\rightarrow 
C^\infty (\Lambda^2M)$ and its pointwise norm is calculated accordingly in 
an orthonormal frame. Hence the Pontryagin number of our ALF spaces is finite.

For any real number $0<R<+\infty$ let $\overline{M}_R\subset M$ 
be the truncated manifold with boundary containing all the points of 
$K\subset M$ as well as those $x\in W\subset M$ for which $r(x)\leqq R$. 

\begin{definition}
Let $(M,g)$ be an ALF four-manifold. Take an arbitrary finite
energy {\rm SU$(2)$}-connection $\nabla_A$ on a (necessarily trivial) rank 
$2$ complex {\rm SU$(2)$} vector bundle $E_0$ over $M$. This connection is 
said to be {\em admissible} if it satisfies two conditions 
{\rm (\cite[Definitions 2.1 and 2.2]{ete-jar})}:
\begin{itemize}
\item[{\rm (i)}] The first is called the {\em weak holonomy condition} 
and says that to $\nabla_A$ there exist constants $0<R<+\infty$ and 
$0<c(g)<+\infty$, this latter being independent of $R$, and a smooth flat 
{\rm SU$(2)$}-connection $\nabla_\Gamma\vert_W$ on $E_0\vert_W$ along the 
end $W\subset M$ such that there exists a gauge on 
$M\setminus\overline{M}_R\subset W$ satisfying
\[\Vert A-\Gamma\Vert_{L^2_{1,\Gamma}\left( M\setminus\overline{M}_R\right)}
\leqq c\Vert F_A\Vert_{L^2\left( M\setminus\overline{M}_R\right)};\]

\item[{\rm (ii)}] The second condition requires $\nabla_A$ to {\em decay 
rapidly} at infinity i.e.,
\[\lim\limits_{R\rightarrow +\infty}\sqrt{R}\:\Vert 
F_A\Vert_{L^2\left( M\setminus\overline{M}_R\right)}=0.\]
\end{itemize}
\label{feltetelek}
\end{definition}

\noindent Note that admissibility could have been defined for an arbitrary Lie 
group $G$. However for $G={\rm SU}(2)$ this definition is natural and 
admissible self-dual connections are the ``good objects'' to consider in the 
ALF scenario as we argued in \cite{ete-jar} because:
\begin{itemize}

\item[(i)] We demonstrated in \cite[Theorem 2.3]{ete-jar} that if 
$N$ in (\ref{fibralas}) is an arbitrary circle bundle over 
$B_{+\infty}\not\cong S^2, \R P^2$ or a trivial circle bundle over 
$B_{+\infty}\cong S^2, \R P^2$ then for any finite energy connection the 
weak holonomy condition is satisfied. Hence in spite of its analytical shape 
it is in fact a mild topological condition only and is essentially always 
valid (except for instance in the important case of the multi-Taub--NUT 
geometries). On the other hand the rapid decay is indeed a non-trivial 
analytical condition and is somewhat stronger than assuming simply finite 
energy;

\item[(ii)] The energy of any admissible ${\rm SU}(2)$
Yang--Mills instanton $\nabla_A$ belongs to a discrete set 
characterized by Chern--Simons invariants $\tau_N(\Gamma_{+\infty} )$ of 
the infinitely distant boundary $N$. More precisely we know 
(\cite[Theorem 2.2]{ete-jar} or Theorem \ref{regitetel} here) that
\begin{equation}
e=\frac{1}{8\pi^2}\Vert F_A\Vert^2_{L^2(M)}\equiv\tau_N(\Gamma_{+\infty} )
\:\:\:\:\:\mbox{{\rm mod} $\Z$}
\label{energia}
\end{equation}
where $\nabla_\Gamma\vert_W =\dd +\Gamma$ is a flat connection 
(associated to $\nabla_A$ by the weak holonomy condition (i) of 
Definition \ref{feltetelek}) in a smooth gauge in which the limit of its 
restriction $\Gamma_{+\infty}:=\lim\limits_{r\rightarrow 
+\infty}\Gamma\vert_{N\times\{ r\}}$ exists (such gauge indeed exists, 
cf. \cite[Section 2]{ete-jar});

\item[(iii)] The framed moduli spaces $\cm (e,\Gamma )$ of irreducible 
admissible ${\rm SU}(2)$ Yang--Mills instantons over a Ricci-flat ALF space are 
smooth possibly empty manifolds and \cite[Theorem 3.2]{ete-jar}
\begin{equation}
\dim\cm (e,\Gamma )=8\left( e+\tau_N(\Theta_{+\infty} )-
\tau_N(\Gamma_{+\infty} )\right) -3b^-(X)
\label{dimenzio}
\end{equation}
where $\nabla_\Theta\vert_W =\dd +\Theta$ is the trivial flat
connection in a gauge in which both the limit 
$\Theta_{+\infty}:=\lim\limits_{r\rightarrow +\infty}
\Theta\vert_{N\times\{ r\}}$ and $\Gamma_{+\infty}$ exists moreover $X$ 
is the Hausel--Hunsicker--Mazzeo compactification (\cite{hau-hun-maz}, 
also defined here soon) of $(M,g)$ with its 
induced orientation.\footnote{Recall that the {\it framed moduli space}
$\cm (e,\Gamma )$ consists of pairs $([\nabla_A],\Gamma )$ where  
$[\nabla_A]$'s are the $L^2_{2,\Gamma}$ gauge equivalence classes of 
irreducible admissible ${\rm SU}(2)$ self-dual connections of energy 
$e\geqq 0$ and asymptotics given by $\nabla_\Gamma\vert_W$ while $\Gamma$ is a 
fixed smooth gauge at infinity. Forgetting about the fixed framing $\Gamma$ at 
infinity we obtain from $\cm (e, \Gamma )$ the {\it unframed moduli space}
$\widehat{\cm}(e,\Gamma )$ consisting of the gauge equivalence classes 
$[\nabla_A]$ only. We therefore find $\dim\widehat{\cm}(e,\Gamma )=
\dim\cm (e, \Gamma )-\dim{\rm SU}(2)=\dim\cm (e, \Gamma )-3$.} 

\end{itemize}

\begin{remark}\rm 1. Several explicit examples demonstrate the relevance of 
both conditions in Definition \ref{feltetelek}. For instance dropping 
only the weak holonomy condition (i) of Definition \ref{feltetelek} there 
exist rapidly decaying smooth reducible ${\rm SU}(2)$ Yang--Mills 
anti-instantons over the multi-Taub--NUT spaces with arbitrary positive energy 
\cite[Sections 2 and 4]{ete-jar} (these will be reviewed in the 
Remark of Sect. \ref{four} here). Likewise, dropping only the rapid 
decay condition (ii) of Definition \ref{feltetelek} there exist smooth 
irreducible ${\rm SU}(2)$ Yang--Mills instantons over the Riemannian 
Schwarzschild space \cite{mos-tav} or over $\R^3\times S^1$ with its flat 
metric \cite{der, nye} with continuous energy spectrum. All of these 
solutions are pathological in some sense: for instance they do not form nice 
moduli spaces.

2. On the contrary, the simplest {\it unframed} moduli space 
$\widehat{\cm}(1,\Theta )$ over the multi-Taub--NUT space consisting of 
all admissible ${\rm SU}(2)$ Yang--Mills anti-instantons with $e=1$ and 
having trivial holonomy at infinity i.e., satisfying $\nabla_\Gamma\vert_W 
=\nabla_\Theta\vert_W$ is a usual moduli space: by (\ref{dimenzio}) it is five 
dimensional and looks like a singular disk fibration over $\R^3$ with usual 
conical singularities corresponding to the reducible solutions at the NUTs. 
It has been constructed explicitly by the aid of the classical conformal 
rescaling method in \cite[Theorem 4.2]{ete-sza}. Moreover all higher integer 
energy moduli spaces $\cm (k,\Theta )$ hence 
$\widehat{\cm}(k,\Theta )$ are non-empty \cite[Theorem 4.3]{ete-sza}.
\end{remark}

\noindent Before proceeding further we take the opportunity and explain
why Definition \ref{feltetelek} is a natural one to impose in the ALF context. 
In this way also the two technical tools which make the ALF scenario so 
special among the non-compact geometries and will be used throughout the 
paper can be introduced. 

The first tool is the {\it Hausel--Hunsicker--Mazzeo compactification} of 
an ALF space \cite{hau-hun-maz}. Take an ALF space $(M,g)$ as 
before. Compactify $M$ by simply shrinking all the circle fibers in the 
fibration (\ref{fibralas}). Let us denote this space by $X$. It is easy 
to see that $X$ is a connected compact smooth $4$-manifold without boundary 
and inherits an orientation from $(M,g)$. Topologically 
\begin{equation}
X=M\cup B_{+\infty}
\label{vegtelen}
\end{equation}
and $B_{+\infty}$ represents a smoothly embedded two codimensional surface 
in $X$. 

The second tool is a {\it codimension $2$ singularity removal theorem of 
Sibner--Sibner} \cite{sib-sib} and {\it R\aa de} \cite{rad}. Let us fix 
some notation which is in agreement with that of \cite{ete-jar}. Take the 
truncated manifold $\overline{M}_R\subset M$ already used in Definition 
\ref{feltetelek} and put $V^\times _R:=M\setminus\overline{M}_R\subseteqq W$ 
for the remaining open tail of the original space $M$. Note that 
$V^\times_R\cong N \times (R,+\infty )$. Consider the fibration 
(\ref{fibralas}) of the boundary and let $U\subset B_{+\infty}$ be a 
coordinate patch of the base space. We obtain a 
corresponding domain $U^\times _R\subset V^\times _R$ what we call an
{\it elementary neighbourhood} \cite{ete-jar}. It follows that
$\pi^{-1}(U)\cong B^2\times S^1$ consequently
$U^\times _R\cong B^2\times S^1\times (R,+\infty )\cong B^2\times 
(B^2)^\times$ i.e., it is a semi-infinite-cylinder-bundle over a disk 
$B^2$ hence $\pi_1(U^\times _R)\cong\Z$. Assume now that a smooth flat local 
${\rm SU}(2)$-connection $\nabla_{\Gamma_m}\vert_{U^\times _R}$ is given.
There exists a canonical gauge $\nabla_{\Gamma_m}\vert_{U^\times _R}=\dd 
+\Gamma_m$ where
\begin{equation}
\Gamma_m=\begin{pmatrix}{\bf i}m & 0\cr
                      0 & -{\bf i}m
          \end{pmatrix}\dd\tau .
\label{lokalismertek}
\end{equation}
Here $\tau\in [ 0,2\pi )$ parameterizes a circle in $U^\times _R$
generating $\pi_1(U^\times _R)\cong\Z$ and $m\in [0,1)$ is the
{\it local holonomy} of the flat local connection. The restriction of the 
globally defined trivial connection $\nabla_\Theta\vert_W$ to $U^\times _R$ 
corresponds to $m=0$ but there may exist further global flat connections 
on $W$ with this property. However in general the flat local connection 
$\nabla_{\Gamma_m}\vert_{U^\times_R}$ does {\it not} extend to a flat 
global connection $\nabla_\Gamma\vert_W$ on the neck.
                      
A tubular neighbourhood $B_{+\infty}\subset V_R\subset X$ of $B_{+\infty}$ in 
$X$ is a $B^2$-bundle over $B_{+\infty}$ and looks like 
$V_R\cong N\times (R,+\infty ]/\sim$ where $\sim$ means that 
$N\times\{+\infty\}$ is pinched into $B_{+\infty}$. Consider a finite open 
covering $B_{+\infty}=\cup_\alpha U_\alpha$ of the base space in 
(\ref{fibralas}) and take the associated elementary neighbourhoods 
$U_{R,\alpha}$ whose collection with $U_{R,\alpha}\cong 
B^2\times S^1\times (R,+\infty ]/\!\!\sim\:\:\cong B^2\times B^2$ gives 
a finite covering for $V_R$. If a finite energy ${\rm SU}(2)$-connection 
$\nabla_A$ is given on $(M,g)$ then the rapid decay condition in 
Definition \ref{feltetelek} makes sure that $\vert F_A\vert_{g'}
\rightarrow 0$ a.e. as $R\rightarrow +\infty$ i.e., it will be a finite energy 
connection on $(X\setminus B_{+\infty},g')$ as well where $g'$ is a 
regularized metric on $X$ which belongs to the conformal class of $g$ on the 
complement of $V_R$ (the original metric does not extend to $X$, even 
conformally). Then the singularity removal theorem we recall now ensures us 
that for sufficiently large $0<R<+\infty$ there exist constants $m\in [0,1)$ 
independent of $\alpha$ and $0<c(g', \alpha )<+\infty$ 
as well as a local $L^2_{1,\Gamma_m}$ gauge on $U^\times_{R,\alpha}$ in 
which $\nabla_A\vert_{U^\times_{R,\alpha}}=\dd +A\vert_{U^\times_{R,\alpha}}$ 
and $\nabla_{\Gamma_m}\vert_{U^\times_{R,\alpha}}=\dd +\Gamma_m$ 
such that
\[\Vert A\vert_{U^\times_{R,\alpha}}-\Gamma_m\Vert_{L^2_{1,\Gamma_m}
(U^\times_{R,\alpha})}\leqq c(g',\alpha )\Vert 
F_A\Vert_{L^2(U_{R,\alpha})}.\]
We recognize this as the local version of the weak holonomy condition
in Definition \ref{feltetelek}. The Sibner--Sibner--R\aa de theorem also says 
that $\nabla_A$ extends over the singularity $B_{+\infty}$ if and only if 
$m=0$. Assume now that the flat local connections 
$\nabla_{\Gamma_m}\vert_{U^\times_{R,\alpha}}$ have a (possibly not 
unique) extensions over the whole $V^\times_R$ and patch together into a 
smooth flat connection $\nabla_\Gamma\vert_{V^\times_R}$. It easily 
follows that the only obstruction against this is the situation if 
$\nabla_{\Gamma_m}\vert_{U^\times_{R,\alpha}}$ with $m\not=0$ is in the 
kernel of $i_*:\pi_1(U^\times_{R,\alpha})\rightarrow \pi_1(V^\times_R)$ 
induced by $i:U^\times_{R,\alpha}\subset V^\times_R$. In 
\cite[Theorem 2.3]{ete-jar} this is converted into an 
easily decidable mild topological condition on the fibration (\ref{fibralas}). 
In other words flat local connections essentially always can be extended 
over the whole neck. Then there exists a gauge on $E_0\vert_{V^\times_R}$ in 
which $\nabla_A\vert_{V^\times_R}=\dd +A$ with $A$ at least in $L^2_{1,\Gamma}$ 
and $\nabla_\Gamma\vert_{V^\times_R}=\dd +\Gamma$ with $\Gamma$ being smooth 
such that one can find smooth gauge transformations independent of $R$ 
satisfying $\gamma^{-1}_\alpha\Gamma_m\gamma_\alpha +
\gamma^{-1}_\alpha\dd\gamma_\alpha=\Gamma\vert_{U^\times _{R,\alpha}}$
where $\Gamma_m$ is the canonical local gauge (\ref{lokalismertek}). 
Then one easily concludes that the local estimates above patch together and 
give part (i) of Definition \ref{feltetelek}. We are convinced now that the 
admissibility condition is essentially a consequence of finite energy.


\section{An improved energy identity for ${\rm SU}(2)$}
\label{three}


After getting some feeling of the admissibility assumption in the case 
of $G={\rm SU}(2)$ Yang--Mills theory over an ALF space, in this section 
we demonstrate that there exist no admissible ${\rm SU}(2)$ 
Yang--Mills instantons of fractional energy over any ALF space i.e., the 
energy identity (\ref{energia}) is superfluous. For the sake of 
completeness first we reproduce \cite[Theorem 2.2]{ete-jar} here. 

Taking into account that the self-duality equations are conformally invariant, 
we can rescale our metric without affecting self-duality. Hence rescale the 
original ALF metric $g$ with a positive function $f:M\rightarrow\R^+$ 
satisfying $f\vert_W\sim O(r^{-2})$ and write $\tilde{g}:=f^2g$. In what 
follows {\it this} rescaled metric $\tilde{g}$ will be used everywhere to 
calculate various Sobolev norms. 

Using the notation of Sect. \ref{two} let $\overline{M}_R\subset M$ be the 
truncated manifold-with-boundary with $r(x)\leqq R$ whose boundary is 
$\partial\overline{M}_R\cong N\times\{ R\}$. Given an 
${\rm SU}(2)$-connection $\nabla_A=\dd +A$ in some    
gauge on the trivial bundle $E_0$ the Chern--Simons functional 
evaluated on its restriction to $\partial\overline{M}_R$ is
\[\tau_{\partial\overline{M}_R}(A_R):=-\frac{1}{8\pi^2}\int
\limits_{\partial\overline{M}_R}{\rm tr}\left(\dd A_R\wedge A_R+
\frac{2}{3}A_R\wedge A_R\wedge A_R\right)\]
in the induced gauge $A_R:=A\vert_{\partial\overline{M}_R}$ on the boundary.

First---motivated by \cite{weh}---we prove two continuity results 
for the Chern--Simons functional in three dimensions which 
are interesting on their own right (also cf. \cite[Lemma 2.1]{ete-jar}). 
Only the first one will be used in this paper.
 
\begin{lemma} Take two ${\rm SU}(2)$-connections $\nabla_{A_R}$ and 
$\nabla_{B_R}$ on the trivial bundle $E_0\vert_{\partial\overline{M}_R}$. 
For some fixed $0<R<+\infty$ on the compact Riemannian three-manifold 
$(\partial\overline{M}_R, \tilde{g}\vert_{\partial\overline{M}_R})$ consider 
Sobolev norms $\Vert\cdot\Vert_{L^p_k}$ with respect to the metric 
$\tilde{g}\vert_{\partial\overline{M}_R}$ and the 
connection $\nabla_{B_R}$ for instance. 

{\rm (i)} Assume that there exists a gauge in which $\nabla_{A_R}=\dd 
+A_R$ and $\nabla_{B_R}=\dd +B_R$ satisfy $A_R, B_R\in 
L^2_1(\partial\overline{M}_R\:;\:   
\Lambda^1(\partial\overline{M}_R)\otimes{\mathfrak s}{\mathfrak u}(2))$. 
Then there exists an estimate
\[\left\vert \tau_{\partial\overline{M}_R}(A_R)\!-\!  
\tau_{\partial\overline{M}_R}(B_R)\right\vert\!\!\leqq\!\!\left(\Vert 
F_{A_R}\Vert_{L^2(\partial\overline{M}_R)}\!+\!
\Vert F_{B_R}\Vert_{L^2(\partial\overline{M}_R)}\right)\!\Vert A_R 
-B_R\Vert_{L^2_1(\partial\overline{M}_R)}+c_1^3\Vert A_R 
-B_R\Vert^3_{L^2_1(\partial\overline{M}_R)}\]
i.e., the Chern--Simons functional is continuous in the $L^2_{1,B_R}$-norm 
in this sense. Here the constant $0<c_1(B_R,R)<+\infty$ is the constant of 
the Sobolev embedding $L^2_1\subset L^3$.

{\rm (ii)} Assume that there exists a gauge in which $\nabla_{A_R}=\dd 
+A_R$ and $\nabla_{B_R}=\dd +B_R$ satisfy $A_R, B_R\in 
L^{\frac{3}{2}}_1(\partial\overline{M}_R\:;\: 
\Lambda^1(\partial\overline{M}_R)\otimes{\mathfrak s}{\mathfrak u}(2))$.
Then there exists a sharper estimate
\[\left\vert \tau_{\partial\overline{M}_R}(A_R)-
\tau_{\partial\overline{M}_R}(B_R)\right\vert\leqq 
c_2\Vert A_R -B_R\Vert_{L^{\frac{3}{2}}_1(\partial\overline{M}_R)}+ 
c_3\Vert A_R -B_R\Vert^3_{L^{\frac{3}{2}}_1(\partial\overline{M}_R)}\]
i.e., the Chern--Simons functional is continuous even in the stronger 
$L^{\frac{3}{2}}_{1,B_R}$-norm in this sense. Here the constants are
\begin{eqnarray}
c_2(B_R,R)&:=&c_4+( c_4+c_4^3)\:\Vert 
B_R\Vert_{L^{\frac{3}{2}}_1(\partial\overline{M}_R)}
+c^3_4\:\Vert B_R\Vert^2_{L^{\frac{3}{2}}_1(\partial\overline{M}_R)}\nonumber\\
&&\nonumber\\
c_3(B_R,R)&:=&c_4+c_4^3+c_4^3\:\Vert 
B_R\Vert_{L^{\frac{3}{2}}_1(\partial\overline{M}_R)}
\nonumber
\end{eqnarray}
with $0<c_4(B_R, R)<+\infty$ being the constant of the sharp Sobolev 
embedding $L^{\frac{3}{2}}_1\subset L^3$.
\label{cslemma}
\end{lemma}

\noindent{\it Proof.} Both inequalities rest on the identity
\[\tau_{\partial\overline{M}_R}(A_R )-\tau_{\partial\overline{M}_R}(B_R )=\]
\begin{equation}
-\frac{1}{8\pi^2}\int\limits_{\partial\overline{M}_R}{\rm tr}
\left( (F_{A_R}+F_{B_R})\wedge (A_R -B_R)-\frac{1}{3}(A_R -B_R )
\wedge (A_R -B_R )\wedge (A_R -B_R )\right) .
\label{azonossag}
\end{equation}

(i) By the aid of H\"older's inequalities with $1=\frac{1}{2}+\frac{1}{2}$ and 
$1=\frac{1}{3}+\frac{1}{3}+\frac{1}{3}$ and the Sobolev embedding 
$L^2_1\subset L^3$ (valid in three dimensions) with $0<c_1(B_R,R)<+\infty$ 
it easily follows that the absolute value of the quadratic term can be 
estimated from above simply by
\begin{eqnarray}
&&\left(\Vert F_{A_R}\Vert_{L^2(\partial\overline{M}_R)}\!+\!
\Vert F_{B_R}\Vert_{L^2(\partial\overline{M}_R)}\right)\!\Vert A_R
-B_R\Vert_{L^2(\partial\overline{M}_R)}\nonumber\\
&&\nonumber\\
&\leqq&\left(\Vert F_{A_R}\Vert_{L^2(\partial\overline{M}_R)}\!+\!
\Vert F_{B_R}\Vert_{L^2(\partial\overline{M}_R)}\right)\!\Vert A_R
-B_R\Vert_{L^2_1(\partial\overline{M}_R)}\nonumber
\end{eqnarray}
while the absolute value of the cubic term has an estimate from above like
\[\Vert A_R -B_R\Vert^3_{L^3(\partial\overline{M}_R)}\leqq c_1^3\Vert A_R 
-B_R\Vert^3_{L^2_1(\partial\overline{M}_R)}\]
which proves the first part.

(ii) Inserting $F_{A_R}=\dd A_R+A_R\wedge A_R$ and
$F_{B_R}=\dd B_R+B_R\wedge B_R$ into (\ref{azonossag}) we 
write $\tau_{\partial\overline{M}_R}(A_R )  
-\tau_{\partial\overline{M}_R}(B_R )$ as the sum of the following three terms:
\[-\frac{1}{8\pi^2}\int\limits_{\partial\overline{M}_R}{\rm tr}
\left(\dd (A_R-B_R)\wedge (A_R-B_R)\right)\]
and 
\[-\frac{1}{8\pi^2}\cdot\frac{2}{3}\int\limits_{\partial\overline{M}_R}
{\rm tr}\left( (A_R -B_R )\wedge (A_R -B_R )\wedge (A_R -B_R )\right)\]
and
\[-\frac{1}{8\pi^2}\int\limits_{\partial\overline{M}_R}{\rm tr}
\left( (2\:\dd B_R+A_R\wedge B_R+B_R\wedge A_R)\wedge (A_R-B_R)\right) .\]
Making use of H\"older's inequalities with $1=\frac{2}{3}+\frac{1}{3}$ 
and $1=\frac{1}{3}+\frac{1}{3}+\frac{1}{3}$, the sharp Sobolev embedding 
$L^{\frac{3}{2}}_1\subset L^3$ (sharply valid in three dimensions) with 
a constant $0<c_4(B_R,R)<+\infty$ and the elementary inequality 
$x^2\leqq x+x^3$ we proceed as follows. The absolute value of the first term 
can be estimated from above by
\begin{eqnarray}
\Vert\dd (A_R-B_R)\Vert_{L^{\frac{3}{2}}(\partial\overline{M}_R)}
\Vert A_R-B_R\Vert_{L^3(\partial\overline{M}_R)}&\leqq &
c_4\Vert A_R-B_R\Vert_{L^{\frac{3}{2}}_1(\partial\overline{M}_R)}
\Vert A_R-B_R\Vert_{L^{\frac{3}{2}}_1(\partial\overline{M}_R)}\nonumber\\
&\leqq&c_4\left(\Vert A_R-B_R\Vert_{L^{\frac{3}{2}}_1(\partial\overline{M}_R)} 
+\Vert A_R-B_R\Vert^3_{L^{\frac{3}{2}}_1(\partial\overline{M}_R)}\right) .
\nonumber
\end{eqnarray}
The absolute value of the second term can be estimated from above by
\[\Vert A_R-B_R\Vert^3_{L^3(\partial\overline{M}_R)}\leqq
c_4^3\Vert A_R-B_R\Vert^3_{L^{\frac{3}{2}}_1(\partial\overline{M}_R)}.\]
Finally we adjust the third term via Stokes' theorem and some algebra 
into the shape
\[-\frac{1}{8\pi^2}\cdot 2\int\limits_{\partial\overline{M}_R}{\rm tr}
\left( B_R\wedge\dd (A_R-B_R)+B_R\wedge (A_R-B_R)\wedge (A_R-B_R)+
B_R\wedge B_R\wedge (A_R-B_R)\right) .\]
Then we can estimate the absolute value of the third term from above by
\begin{eqnarray}
&&\Vert B_R\Vert_{L^3(\partial\overline{M}_R)}\Vert \dd 
(A_R-B_R)\Vert_{L^{\frac{3}{2}}(\partial\overline{M}_R)}+
\Vert B_R\Vert_{L^3(\partial\overline{M}_R)}
\Vert A_R-B_R\Vert^2_{L^3(\partial\overline{M}_R)}\nonumber\\
&&\nonumber\\
&&+\Vert B_R\Vert^2_{L^3(\partial\overline{M}_R)}                   
\Vert A_R-B_R\Vert_{L^3(\partial\overline{M}_R)}\nonumber\\
&&\nonumber\\
&\leqq&\left( (c_4+c^3_4)\Vert B_R\Vert_{L^{\frac{3}{2}}_1
(\partial\overline{M}_R)}+c^3_4\Vert 
B_R\Vert^2_{L^{\frac{3}{2}}_1(\partial\overline{M}_R)}\right)
\Vert(A_R-B_R)\Vert_{L^{\frac{3}{2}}_1(\partial\overline{M}_R)}\nonumber\\
&&\nonumber\\
&&+c_4^3\Vert B_R\Vert_{L^{\frac{3}{2}}_1(\partial\overline{M}_R)}
\Vert(A_R-B_R)\Vert^3_{L^{\frac{3}{2}}_1(\partial\overline{M}_R)}.\nonumber
\end{eqnarray}
Putting together these estimates we obtain the inequality of the 
second part of the lemma. 
$\Diamond$

\begin{remark}\rm We also record here that in the case of flat 
connections (\ref{azonossag}) shows that the inequalities of the lemma 
cut down to
\[\left\vert\tau_{\partial\overline{M}_R}(\Gamma_R')-\tau_{\partial
\overline{M}_R}(\Gamma_R'')\right\vert\leqq
\Vert\Gamma_R'-\Gamma_R''\Vert^3_{L^3(\partial\overline{M}_R)}\leqq
c_1^3\Vert\Gamma_R'-\Gamma_R''\Vert^3_{L^2_1 (\partial\overline{M}_R)}\]
and a similar one for the $L^{\frac{3}{2}}_1$-norm.
\end{remark}

\noindent Now we are in a position to reprove \cite[Theorem 2.2]{ete-jar}
following the steps of \cite[Section 2]{ete-jar}. Take an admissible ${\rm 
SU}(2)$-connection $\nabla_A$ on $E_0$ and the 
corresponding flat connection $\nabla_\Gamma\vert_{V^\times _R}$ to which it 
converges. Suppose that we are in the gauge on $E_0\vert_{V^\times _R}$ in 
which $\nabla_A\vert_{V^\times _R}=\dd +A$ and 
$\nabla_\Gamma\vert_{V^\times _R}=\dd +\Gamma$ and the corresponding 
connection $1$-forms satisfy the inequality in part (i) of Definition 
\ref{feltetelek}. Let $A_R$ and $\Gamma_R$ be their restrictions to 
$\partial\overline{M}_R$. 

\begin{theorem} {\rm (cf. \cite[Theorem 2.2]{ete-jar})} Let $(M,g)$ be 
an ALF space with an end $W\cong N\times\R^+$. Let $E_0$ be an ${\rm SU}(2)$ 
vector bundle over $M$, necessarily trivial, with an admissible ${\rm 
SU}(2)$ Yang--Mills instanton $\nabla_A$ on it: a smooth, finite 
energy self-dual connection satisfying Definition \ref{feltetelek}. Then
\[\frac{1}{8\pi^2}\Vert 
F_A\Vert^2_{L^2(M)}\equiv\tau_N(\Gamma_{+\infty})\:\:\:\:\:{\rm mod}\:\Z\]
that is, its energy is congruent to a Chern--Simons invariant of the 
boundary given by the flat connection $\nabla_\Gamma\vert_W$ in part 
(i) of Definition \ref{feltetelek}. 
\label{regitetel}
\end{theorem}

\noindent {\it Proof.} We estimate the difference 
$\left\vert\tau_{\partial\overline{M}_R}(A_R)-
\tau_{\partial\overline{M}_R}(\Gamma_R)\right\vert^{\frac{2}{3}}$ along 
$V^\times_R\cong N\times (R,+\infty )$ as follows. As a first step 
by the aid of the mean value theorem we find an $R_0\in (R, 2R)$ such that 
for all $2R\leqq S<+\infty$
\[\frac{1}{2R}\left\vert\tau_{\partial\overline{M}_{R_0}}(A_{R_0})-
\tau_{\partial\overline{M}_{R_0}}(\Gamma_{R_0})\right\vert^{\frac{2}{3}}\leqq
\int\limits_R^S\left\vert\tau_{\partial\overline{M}_r}(A_r)-
\tau_{\partial\overline{M}_r}(\Gamma_r)\right\vert^{\frac{2}{3}}r^{-2}\dd r\]
holds; then we go on via the inequality of part (i) of Lemma \ref{cslemma} with 
$F_{\Gamma_r}=0$ and the elementary inequality $(x+y)^{\frac{2}{3}}\leqq 
2(x^{\frac{2}{3}}+y^{\frac{2}{3}})$ to get 
\begin{eqnarray}
&\leqq &\int\limits_R^S
\left( \Vert F_{A_r}\Vert_{L^2(\partial\overline{M}_r)}\Vert 
A_r-\Gamma_r\Vert_{L^2_{1,\Gamma_r}
(\partial\overline{M}_r)}+c_1^3\Vert A_r-\Gamma_r\Vert^3_
{L^2_{1,\Gamma_r}(\partial\overline{M}_r)}\right)^{\frac{2}{3}}r^{-2}\dd r
\nonumber\\
&\leqq &2\int\limits_R^S\left(\Vert 
F_{A_r}\Vert^{\frac{2}{3}}_{L^2(\partial\overline{M}_r)}\Vert 
A_r-\Gamma_r\Vert^{\frac{2}{3}}_{L^2_{1,\Gamma_r}
(\partial\overline{M}_r)}+c_1^2\Vert A_r-\Gamma_r\Vert^2_{L^2_{1,\Gamma_r}
(\partial\overline{M}_r)}\right)r^{-2}\dd r\:;\nonumber
\end{eqnarray}
and then make further steps by applying on the first term two H\"older's 
inequalities along $((R,S),r^{-2}\dd r )$ with $1=\frac{2}{3}+\frac{1}{3}$ 
and then with $1=\frac{1}{2}+\frac{1}{2}$ to obtain 
\begin{eqnarray}
&\leqq &2c_5^{\frac{1}{3}}
\Vert F_A\Vert^{\frac{4}{3}}_{L^2(V^\times_R\setminus V^\times_S)}
\left(\:\int\limits_R^S\Vert A_r-\Gamma_r\Vert^2_{L^2_{1,\Gamma_r}
(\partial\overline{M}_r)}r^{-2}\dd r\right)^{\frac{1}{3}}
+2\!\!\int\limits_R^Sc_1^2\Vert A_r-
\Gamma_r\Vert^2_{L^2_{1,\Gamma_r}(\partial\overline{M}_r)}r^{-2}\dd r\nonumber\\
&\leqq &2c_5\Vert F_A\Vert^{\frac{4}{3}}_{L^2(V^\times_R\setminus V^\times_S)}
\Vert A-\Gamma\Vert^{\frac{2}{3}}_{L^2_{1,\Gamma}
(V^\times _R\setminus V^\times_S )}+2c_6\Vert A-
\Gamma\Vert^2_{L^2_{1,\Gamma}(V^\times_R\setminus V^\times_S)}\nonumber\\
&\leqq& 2c_5\Vert F_A\Vert^{\frac{4}{3}}_{L^2(V^\times_R)}
\Vert A-\Gamma\Vert^{\frac{2}{3}}_{L^2_{1,\Gamma}(V^\times _R)}+2c_6\Vert A-   
\Gamma\Vert^2_{L^2_{1,\Gamma}(V^\times_R)}.\nonumber
\end{eqnarray}
Here the constant $0<c_5(\tilde{g})<+\infty$, independent of $R$ takes 
into account the discrepancy between the original measure on 
$V^\times_R$ and its restriction to $\partial\overline{M}_r$ multiplied by 
$r^{-2}\dd r$ (cf. the asymptotical shape of an ALF metric in Sect. \ref{two}). 
Moreoever $c_6:=c_5\sup_{r\in [R,+\infty ]}c_1^2(r)<+\infty$. 

Putting these steps together and then referring to part (i) 
of Definition \ref{feltetelek} we therefore come up with the estimate
\[\left\vert\tau_{\partial\overline{M}_{R_0}}(A_{R_0})-
\tau_{\partial\overline{M}_{R_0}}(\Gamma_{R_0})\right\vert
\leqq (8c_7)^{\frac{3}{2}}\left(R\Vert
F_A\Vert^2_{L^2(V^\times_R)}\right)^{\frac{3}{2}}=
(8c_7)^{\frac{3}{2}}\left(\sqrt{R}\:\Vert
F_A\Vert_{L^2(V^\times_R)}\right)^3\]
with some $R_0\in (R,2R)$. Checking the shape of the constants in part (i) 
of Definition \ref{feltetelek} 
and part (i) of Lemma \ref{cslemma} we can assume that the overall 
constant $0<c_7(\tilde{g},\Gamma )<+\infty$ is bounded. Hence we obtain from 
this last inequality by the aid of the rapid decay 
condition i.e., part (ii) of Definition \ref{feltetelek} that
\[\lim\limits_{R\rightarrow +\infty}\vert\tau_{\partial\overline{M}_R}(A_R)-
\tau_{\partial\overline{M}_R}(\Gamma_R)\vert =0\]
or writing $\Gamma_{+\infty}=\lim\limits_{R\rightarrow 
+\infty}\Gamma_R$ and regarding this as a flat connection on 
the infinitely distant boundary $N$ of (\ref{fibralas}) we conclude that 
\[\lim\limits_{R\rightarrow +\infty}\tau_{\partial\overline{M}_R}(A_R) 
=\tau_N(\Gamma_{+\infty})\]
which gives the result when applied to a self-dual admissible ${\rm SU}(2)$ 
connection $\nabla_A$ over the ALF space $(M,g)$ as claimed. $\Diamond$
\vspace{0.1in}

\noindent It follows from \cite[Theorem 4.3]{kir-kla} already at this 
point that the energy spectrum consists of rational numbers. Now we show that
in fact only those flat connections appear on which the Chern--Simons 
functional takes integer values i.e., $\tau_N(\Gamma_{+\infty})\in\Z$ in the 
previous theorem.

\begin{theorem}
Let $(M,g)$ be an ALF space with an end $W\cong N\times\R^+$. Let $E_0$ be
an ${\rm SU}(2)$ vector bundle over $M$, necessarily trivial, with 
an admissible ${\rm SU}(2)$ Yang--Mills instanton $\nabla_A$ on it i.e., a 
smooth, finite energy self-dual connection satisfying Definition 
\ref{feltetelek}. Then 
\[\frac{1}{8\pi^2}\Vert F_A\Vert^2_{L^2(M)}\in\N\]
that is, in addition to Theorem \ref{regitetel} its energy is always integer. 

Regarding the asymptotical shape of $\nabla_A$ if $M$ is in addition 
simply connected then the associated flat connection $\nabla_\Gamma\vert_W$ 
of part (i) in Definition \ref{feltetelek} has trivial local holonomy at 
infinity i.e., $m=0$ in (\ref{lokalismertek}) 
(in this case if $\nabla_\Gamma\vert_W\not=\nabla_\Theta\vert_W$ then 
$\pi_1(B_{+\infty})\not=1$ in (\ref{fibralas}) and (\ref{vegtelen})). 
\label{fotetel}
\end{theorem}

\noindent{\it Proof.} We continue to estimate $\tau_N(\Gamma_{+\infty})$ in 
Theorem \ref{regitetel} by exploiting the consequences of the rapid decay 
condition i.e., part (ii) of Definition \ref{feltetelek} more systematically. 
The key idea is to work over the Hausel--Hunsicker--Mazzeo compactification 
(\ref{vegtelen}). 

A tubular neighbourhood $B_{+\infty}\subset V_R\subset X$ of the 
infinite distant surface then looks like $V_R =V^\times_R\cup B_{+\infty}$. We 
know by regularity of Yang--Mills fields that for sufficiently large $R$ we 
can suppose that the gauge we use is smooth. Applying the mean value theorem 
as in the proof of Theorem \ref{regitetel} we get for some $R_0\in (R,2R)$ 
that 
\[\Vert
A-\Gamma\Vert_{L^2_{1,\Gamma}(\partial\overline{M}_{R_0})}\leqq
\sqrt{2c_5R}\:c\Vert F_A\Vert_{L^2(V^\times_R)},\:\:\:\:\:
\Vert F_A\Vert_{L^2(\partial\overline{M}_{R_0})}\leqq \sqrt{2c_5R}\:
\Vert F_A\Vert_{L^2(V^\times_R)}.\] 
Therefore by part (ii) of Definition \ref{feltetelek} we find that 
$(A-\Gamma )\vert_{B_{+\infty}}=0$ and 
$\nabla_\Gamma (A-\Gamma )\vert_{B_{+\infty}}=0$ as well as 
$F_A\vert_{B_{+\infty}}=0$ hold pointwise. 
Consequently putting $\omega :={\rm tr}(\dd A\wedge A+\frac{2}{3}A\wedge
A\wedge A)$ and referring to the identity $\dd\omega = 
{\rm tr}(F_A\wedge F_A)$ we conclude that 
$\omega\in L^2_{1,\Gamma}(\overline{V}_R\:;\:\Lambda^3\overline{V}_R)$ where 
$\overline{V}_R$ is the closure with 
$\partial\overline{V}_R =(\partial\overline{M}_R)^*$ (reversed orientation). 
Moreover $\omega$ is smooth near $\partial\overline{V}_R$. Hence by (the 
generalized) Stokes' theorem there exists an integer $n(R)\in\Z$ such that
\[\left\vert n(R)-\tau_{\partial\overline{M}_R}(A_R)\right\vert 
=\frac{1}{8\pi^2}\left\vert\:\int\limits_{V_R}{\rm tr}\left( F_A\wedge 
F_A\right)\right\vert\leqq\Vert F_A\Vert^2_{L^2(V^\times_R)}.\]
By Theorem \ref{regitetel} we know that 
$\tau_{\partial\overline{M}_R}(A_R)\rightarrow \tau_N(\Gamma_{+\infty})$ 
as $R\rightarrow +\infty$ and the right hand side of this last equation can 
be kept as small as we please as $R\rightarrow +\infty$ since $\nabla_A$ has 
finite energy. We conclude that $\tau_N(\Gamma_{+\infty})\in\Z$ in 
Theorem \ref{regitetel}.

Consider an elementary neighbourhood $U^\times_R\subset V^\times_R$. If 
$M$ is in addition supposed to be simply connected then any loop representing 
a generator of $\pi_1(U^\times_R)\cong\Z$ is contractible in $M$. In other 
words if $\ell_0:S^1\rightarrow U^\times_R\subset M$ is a non-trivial loop 
then it shrinks to the trivial one $\ell_1(S^1)=x_0\in M$. Take a smooth 
non-self-intersecting curve 
$\gamma : [\frac{R}{2}, +\infty )\rightarrow M\subset X$ connecting 
$\gamma (\frac{R}{2})=x_0$ with the limit point 
$\gamma (+\infty )\in B_{+\infty }$. If $r\in [\frac{R}{2}, +\infty )$ 
denotes the coordinate along this curve then we can suppose that it coincides 
with the radial coordinate we use along 
$U^\times_R\cong B^2\times S^1\times (R,+\infty )$. In an open tubular 
neighbourhood $O$ of $\gamma$ in $X$ we can suppose that we are in a radial 
gauge i.e., $\nabla_A\vert_O=\dd +A\vert_O$ satisfies $A_r=0$. Therefore 
in this gauge $\frac{\partial A_\tau}{\partial r} =F_{r\tau}$ and we also 
have $A_\tau (\gamma (\frac{R}{2}))=A_\tau\vert_{x_0}=0$ because $\nabla_A$ is 
smooth in $x_0\in M$. Consequently integrating $A_\tau$ along 
$[\frac{R}{2},+\infty )$ we obtain an estimate
\[\sqrt{2}\:m=\vert\Gamma_m\vert_{g'}=\vert\Gamma_\tau (\gamma (+\infty ))
\vert_{g'}=\vert A_\tau (\gamma (+\infty ))\vert_{g'}\leqq 
\frac{1}{R}\sup\limits_\rho\vert F_{r\tau}(\rho)\vert_{g'} <+\infty\]
since $F_{r\tau}$ is bounded. But $R$ is arbitrary hence $m=0$ in 
(\ref{lokalismertek}) implying $\nabla_\Gamma\vert_{V^\times_R}$ is a 
flat connection with trivial local holonomy at infinity. 

We conclude that the admissible self-dual connection $\nabla_A$ has integer 
energy. $\Diamond$

\begin{remark}\rm 1. If $\nabla_\Gamma\vert_{V^\times_R}$ 
has trivial local holonomy at infinity then it is the pullback of a flat 
connection via $\pi\times{\rm Id}_{(R,+\infty )}:V^\times_R\cong N\times 
(R,+\infty )\rightarrow B_{+\infty}\times (R,+\infty )$. Consequently the
restriction $\Gamma_{+\infty}$ of $\Gamma\vert_{V^\times_R}$ to the
infinitely distant boundary $N$ in (\ref{fibralas}) is also a pullback
connection $1$-form via $\pi :N\rightarrow B_{+\infty}$. It is known
that such a flat connection has vanishing Chern--Simons invariant
(cf. \cite[pp. 547-548 and Theorem 4.3]{kir-kla}) i.e., we find independently 
of our considerations above that in the special case of simply connected 
$M$ the expression $\tau_N(\Gamma_{+\infty})$ is an integer.

2. We cannot achieve in general that the local holonomy at 
infinity vanishes. For example there exists a family of smooth flat hence 
self-dual admissible connections on $\R^3\times S^1$ parameterized by their 
holonomy $m\in [0,1)$. Hence $m$ is also their holonomy at the infinite 
$N\cong S^2\times S^1$.
\end{remark}



\section{The case of ${\rm U}(2)$}
\label{four}


In this section we demonstrate that Theorem \ref{fotetel} 
continues to hold for the slightly larger gauge groups 
${\rm SU}(2)\times {\rm U}(1)$ and ${\rm U}(2)$ but in the latter case 
with a topological condition on the underlying vector bundle. 

However before doing this we make an important comment here. As we already 
mentioned, admissibility as it stands in Definition \ref{feltetelek} can 
be formulated for an arbitrary (compact) Lie group $G$. Hence repeating 
the proof of the previous section we could obtain an analogue of Theorem 
\ref{fotetel} for arbitrary $G$. The reason we do not do this is 
that for Yang--Mills instantons with a general Lie group the admissibility 
would indeed be a very strong assumption hence the analogue of 
Theorem \ref{fotetel} would be a rather weak statement. This is because for 
general $G$ the analogue of the powerful singularity removal 
theorem of Sibner--Sibner \cite{sib-sib} and R\aa de \cite{rad} is not 
known consequently the weak holonomy condition might in principle be a 
very strong requirement, cf. our discussion in Sect. 2 above.

Rather we restrict ourselves to those Lie groups which can be somehow 
``traced back'' to ${\rm SU}(2)$ in order to keep our results strong.

\begin{lemma}
Let $(M,g)$ be an ALF space and $\tilde{E}$ be a rank $2$ complex 
${\rm SU}(2)\times {\rm U}(1)$ vector bundle on $M$. 

{\rm (i)} Then $\tilde{E}\cong E_0\otimes L$ where $E_0$ is a 
(necessarily trivial) rank 2 complex ${\rm SU}(2)$ vector bundle and $L$ is 
a ${\rm U}(1)$ line bundle over $M$;

{\rm (ii)} Every ${\rm SU}(2)\times {\rm U}(1)$-connection on 
$\tilde{E}\cong E_0\otimes L$ is of the form 
\[\nabla_A\otimes{\rm Id}_L+{\rm Id}_{E_0}\otimes\nabla_B\]
where $\nabla_A$ is an ${\rm SU}(2)$-connection on $E_0$ 
and $\nabla_B$ is an ${\rm U}(1)$-connection on $L$;

{\rm (iii)} The curvature of this product connection looks like
\[F_A+(F_B\oplus F_B)\]
hence an ${\rm SU}(2)\times {\rm U}(1)$-connection on $\tilde{E}$ is 
self-dual if and only if both 
its ${\rm SU}(2)$ and ${\rm U}(1)$ parts are self-dual;

{\rm (iv)} If an ${\rm SU}(2)\times {\rm U}(1)$-connection on $\tilde{E}$ 
has finite energy $e$ then it admits a decomposition
\[e=\frac{1}{8\pi^2}\Vert F_A\Vert^2_{L^2(M)}+\frac{1}{8\pi^2}\Vert 
F_B\oplus F_B\Vert^2_{L^2(M)}=\frac{1}{8\pi^2}\Vert 
F_A\Vert^2_{L^2(M)}+\frac{1}{4\pi^2}\Vert F_B\Vert^2_{L^2(M)}.\]
\label{u2lemma}
\end{lemma}
\noindent {\it Proof.} (i) Standard obstruction theory says that over a 
non-compact oriented four-manifold $M$ principal bundles with a 
connected, compact structure group $G$ are classified by $H^2(M;\pi_1(G))$. 
Hence on the one hand all principal ${\rm SU}(2)$-bundles are trivial 
over $M$ since $\pi_1({\rm SU}(2))\cong 1$. Let us denote the associated 
complex rank $2$ trivial bundle by $E_0$. On the other hand principal 
${\rm U}(1)$-bundles are classified by $H^2(M;\Z )$ since 
$\pi_1({\rm U}(1))\cong\Z$. Referring to the canonical isomorphism 
$\pi_1({\rm U}(1))\cong\pi_1({\rm SU}(2)\times{\rm U}(1))$ we obtain that 
a rank $2$ complex ${\rm SU}(2)\times{\rm U}(1)$ vector bundle $\tilde{E}$ 
can be uniquely written in the form $E_0\otimes L$ where $L$ is a 
${\rm U}(1)$ line bundle. 

(ii) Let $\nabla_A$ be an ${\rm SU}(2)$-connection on $E_0$ and 
$\nabla_B$ be an ${\rm U}(1)$-connection on $L$. Taking the 
embeddings ${\mathfrak s}{\mathfrak u}(2)\subset
{\mathfrak s}{\mathfrak u}(2)\times{\mathfrak u}(1)\cong{\mathfrak u}(2)$ 
as usual and ${\mathfrak u}(1)
\subset{\mathfrak s}{\mathfrak u}(2)\times{\mathfrak u}(1)
\cong{\mathfrak u}(2)$ given by 
$B\mapsto\begin{pmatrix} 
               B & 0\\
               0 & B
         \end{pmatrix}=B\oplus B$ we can identify a product connection on 
$E_0\otimes L$ with an ${\rm SU}(2)\times{\rm U}(1)$-connection 
$\nabla_{A+(B\oplus B)}$ on the corresponding vector bundle $\tilde{E}$ 
and vice versa. 

(iii) Regarding the curvature we calculate
\[F_{A+(B\oplus B)}=F_A+(F_B\oplus F_B)+A\wedge (B\oplus B) 
+(B\oplus B)\wedge A= F_A+(F_B\oplus F_B)\]
since $A\wedge (B\oplus B) +(B\oplus B)\wedge A=0$ because $B\oplus B$ is in 
the centre of ${\mathfrak u}(2)$. Consequently
\[*F_{A+(B\oplus B)}=*(F_A+(F_B\oplus F_B))=*F_A+(*F_B\oplus 
*F_B)=F_A+(F_B\oplus F_B)=F_{A+(B\oplus B)}\]
demonstrating that our product connection continues to be self-dual if 
and only if its components are self-dual. 

(iv) Finally the energy can be decomposed as claimed. $\Diamond$
\vspace{0.1in}

\noindent Consider an ${\rm SU}(2)\times{\rm U}(1)$ Yang--Mills instanton 
over some ALF space $(M,g)$. Then it follows from part (iv) of Lemma 
\ref{u2lemma} that it is admissible\footnote{To be precise for this we 
generalize Definition \ref{feltetelek} i.e., the definition of 
${\rm SU}(2)$-admissibility to arbitrary $G$-admissibility in a 
straightforward way.} if and only if both its ${\rm SU}(2)$ component 
$\nabla_A$ and ${\rm U}(1)$ component $\nabla_B$ are admissible. We already 
know from Theorem \ref{fotetel} that if $\nabla_A$ is admissible then it 
has not only finite but even integer energy: 
\[\frac{1}{8\pi^2}\Vert F_A\Vert^2_{L^2(M)}=k.\]
Regarding the Abelian part we can take the embedding ${\rm 
U}(1)\subset{\rm SU}(2)$ and identify $\nabla_B$ on the line bundle $L$ 
with the reducible ${\rm SU}(2)$ Yang--Mills instanton 
$\nabla_B\mapsto\nabla_B\oplus\nabla_{-B}$ on the split bundle 
$L\oplus L^{-1}$. This is an admissible solution hence we can apply 
either Theorem \ref{regitetel} to conclude that its energy must be 
congruent to an ${\rm U}(1)$ Chern--Simons invariant of the boundary 
(\ref{fibralas}) hence must be integer or we can simply refer to 
Theorem \ref{fotetel} again to obtain that it has integer energy 
$\frac{1}{8\pi^2}\Vert F_B\oplus F_{-B}\Vert^2_{L^2(M)}=l$. Hence
\[\frac{1}{8\pi^2}\Vert F_B\oplus F_B\Vert^2_{L^2(M)}
=\frac{1}{8\pi^2}\Vert F_B\oplus F_{-B}\Vert^2_{L^2(M)}=l.\]
Therefore an admissible ${\rm SU}(2)\times{\rm U}(1)$ instanton 
$\nabla_{A+(B\oplus B)}$ has integer energy by part (iv) of Lemma 
\ref{u2lemma}. 

Finally note that if $E$ is an ${\rm U}(2)$ vector bundle with
$w_2(E^\R )=0\in H^2(M;\Z_2)$ then it uniquely lifts to an
${\rm SU}(2)\times{\rm U}(1)$ vector bundle $\tilde{E}$. Consequently any
connection on $E$ can be uniquely lifted to a connection on $\tilde{E}$.
Taking into account the canonical isomorphism ${\mathfrak u}(2)
\cong{\mathfrak s}{\mathfrak u}(2)\times{\mathfrak u}(1)$ we find that 
the energies of the original and the lifted connections are equal.
Consequently our considerations continue to hold for admissible ${\rm
U}(2)$ Yang--Mills instantons on bundles with vanishing second
Stiefel--Whitney class.

Summing up we obtain

\begin{theorem}
Let $(M,g)$ be an ALF space with an end $W\cong N\times\R^+$. Let $E$ be a 
${\rm U}(2)$ vector bundle with $w_2(E^\R )=0\in H^2(M;\Z_2)$ carrying 
an admissible ${\rm U}(2)$ Yang--Mills instanton $\nabla_A$ i.e., a 
smooth self-dual connection satisfying Definition \ref{feltetelek}. Then
\[\frac{1}{8\pi^2}\Vert F_A\Vert^2_{L^2(M)}\in\N\]
that is, its energy is always integer. 

Regarding the asymptotical shape of $\nabla_A$ if $M$ is in addition
simply connected then the associated flat connection 
$\nabla_\Gamma\vert_W$
of part (i) in Definition \ref{feltetelek} has trivial local holonomy at
infinity i.e., $m=0$ in (\ref{lokalismertek})
(in this case if $\nabla_\Gamma\vert_W\not=\nabla_\Theta\vert_W$ then
$\pi_1(B_{+\infty})\not=1$ in (\ref{fibralas}) and (\ref{vegtelen})). 
$\Diamond$
\label{U2}
\end{theorem}

\begin{remark}\rm We would like to point out the relevance of 
admissibility in the case of Abelian instantons. In this case it is easy 
to see in the framework of $L^2$-cohomology that without imposing 
admissibility a continuous energy spectrum would destroy everything.

So let $\nabla_B$ be an Abelian instanton over an ALF space $(M,g)$. 
By definition its energy is finite but in general 
$[\frac{1}{2\pi{\bf i}}F_B]\in \overline{H}^2_{L^2}(M,g)$ i.e., 
the curvature lives only in the second (reduced) $L^2$-cohomology group of 
the non-compact but complete Riemannian $4$-manifold $(M,g)$ with ALF 
asymptotics. It turns out \cite{hau-hun-maz} that for such geometries 
this subtle cohomology reduces to ordinary cohomology of the 
Hausel--Hunsicker--Mazzeo compactification (\ref{vegtelen}). Therefore 
$\overline{H}^2_{L^2}(M,g)\cong H^2(X;\R )$ and using (\ref{fibralas}) 
and (\ref{vegtelen}) the Mayer--Vietoris sequence gives 
\[\dots\rightarrow H^1(N;\R )\rightarrow 
\overline{H}^2_{L^2}(M,g) \rightarrow H^2_c(M;\R )
\oplus H^2(B_{+\infty};\R )\rightarrow H^2(N;\R )\rightarrow\dots\]
Consequently these cohomology classes can be divided into 
two parts as follows: we say that an $L^2$-cohomology class on a complete 
Riemannian manifold $(M,g)$ is {\it topological} if its image lies in the 
ordinary compactly supported de Rham cohomology $H^2_c(M;\R )$ under the 
homomorphism above. Otherwise it is called {\it non-topological} i.e., if 
its image is in $H^2(B_{+\infty};\R )$. Note that representatives of 
non-topological $L^2$-cohomology classes are necessarily exact $2$-forms on 
$M$. Roughly speaking \cite{seg-sel} non-topological $L^2$-cohomology 
classes are not predictable by topological means. 

We make two assumptions. The first is that $H^j(N;\R )=0$ $(j=1,2)$. In this 
case unambigously 
\[\frac{1}{2\pi{\bf i}}F_B =\sum\limits_{i=1}^{b^2_c(M)}k_i\omega_i +\omega_0\]
where $\{\omega_i\}$ are harmonic representatives of the basis of the compactly 
supported integer cohomology $H^2_c(M;\R )\cap H^2(M;\Z )$ and $k_i\in\Z$ and 
$\omega_0=\dd\beta$ represents the exact non-topological part. Obviously 
$\nabla_B$ lives on a line bundle $L$ with Chern 
class $c_1(L)=[\:\sum k_i\omega_i]\in H^2_c(M;\R )\cap H^2(M;\Z )$. 
The finite energy of this Abelian instanton looks like
\begin{eqnarray}
\frac{1}{8\pi^2}\Vert F_B\Vert^2_{L^2(M)}&=&\frac{1}{2}\int
\limits_M\left(\frac{1}{2\pi{\bf i}}F_B \right)\wedge *
\left(\frac{1}{2\pi{\bf i}}F_B\right)\nonumber\\ 
&=&\frac{1}{2}\sum\limits_{i,j}k_ik_j\int\limits_M\omega_i\wedge *\omega_j+
\sum\limits_ik_i\int\limits_M\omega_0\wedge *\omega_i+
\frac{1}{2}\int\limits_M\omega_0\wedge *\omega_0.\nonumber
\end{eqnarray} 
The second assumption is that all $k_i=0$. In 
this case the curvatue is an exact self-dual non-compactly supported 
$2$-form on $(M,g)$ still having finite energy. The corresponding 
Abelian instanton $\nabla_{B_0}$ with $F_{B_0}=2\pi\ii\omega_0$ lives on the 
trivial bundle $L_0\cong M\times\C$. Taking into account the triviality of 
$L_0$ as well as the Abelian nature of $\nabla_{B_0}$ if it is 
self-dual on $L_0$ with $\frac{1}{8\pi^2}\Vert F_{B_0}\Vert^2_{L^2(M)}=
\frac{1}{2}$ then for any $c\in\R$ the rescaled connection $\nabla_{cB_0}$ 
remains smooth and self-dual on $L_0$ with energy 
$\frac{1}{8\pi^2}\Vert F_{cB_0}\Vert^2_{L^2(M)}=\frac{c^2}{2}$. This 
unexpected continuous energy phenomenon occurs for example in the important 
case of the multi-Taub--NUT spaces \cite[Sections 2 and 4]{ete-jar}: the 
connection $1$-form $B_0$ arises as the metric dual of the Killing field 
associated to the ${\rm U}(1)$ isometry of the metric. 

Consequently in this case admissibility (more precisely part (i) of 
Definition \ref{feltetelek}) is indeed to be imposed 
which gives $c=k$ hence the energy of the ${\rm U}(1)$ Yang--Mills 
instanton $\nabla_{kB_0}$ is a half-integer. Hence the energy of the 
reducible ${\rm SU}(2)$ Yang--Mills instanton $\nabla_{kB_0\oplus 
(-kB_0)}$ or the ${\rm U}(2)$ one $\nabla_{kB_0\oplus kB_0}$ will be 
integer in agreement with Theorems \ref{fotetel} or \ref{U2}.
\end{remark}


\section{Conclusion and outlook toward quantum theory}
\label{five}


In this paper we proved that the energy spectrum of a natural class of 
${\rm SU}(2)$ or ${\rm U}(2)$ Yang--Mills instantons (called admissible 
instantons) over a generic (i.e., not necessarily hyper-K\"ahler) ALF space 
consists of non-negative integers only. This sharpens the previously 
known result that the energy must be congruent to a Chern--Simons invariant of 
the infinitely distant boundary $N$ hence to a rational number. In this 
context the stronger result is surprising because $N$ as defined in 
(\ref{fibralas}) has many non-integer Chern--Simons invariants, cf. 
\cite[Theorem 4.3]{kir-kla}. We have seen that the reason behind this 
integrality phenomenon is the existence of a smooth compactification of the 
original space and a powerful codimension $2$ singularity removal result.

In this closing section we would like to push one step further the role played 
by this Hausel--Hunsicker--Mazzeo compactification and this 
Sibner--Sibner--R\aa de singularity removal theorem in the ALF scenario. 
Namely, we ask ourselves whether the emergence of the infinitely distant 
surface $B_{+\infty}$ in (\ref{fibralas}) and (\ref{vegtelen}) is a 
topological hint that the thing ``lurks'' behind the concept of a four 
dimensional Yang--Mills theory over an ALF space is in fact a $2$ 
dimensional conformal field theory. If this is the case---as we will 
informally argue below---then it would help one to construct the 
underlying (twisted $N=2$ supersymmetric) quantum gauge theory. 

We want to calculate the partition function of our quantum gauge 
theory over the ALF space $(M,g)$. Put
\[\ca (E_0):=\left\{\nabla_A\:\left\vert\:\mbox{$\nabla_A$ 
is admissible on $E_0$ in the sense of Definition 
\ref{feltetelek}}\right\}\right.\]
and let $\cb (E_0):=\ca (E_0)/\cg (E_0)$ be the quotient space of 
$L^2_{2,\Gamma}$ gauge equivalence classes. Taking the complex coupling 
constant $\tau :=\frac{\theta}{2\pi}+\frac{4\pi}{e^2}{\bf i}\in\C^+$ 
and suppressing the supersymmetric terms the partition 
function is a complex number given by the formal integral 
\[Z(M,g,\tau ,{\rm SU}(2))=\int\limits_{\cb (E_0)}{\rm e}^{\frac{1}{2e^2}\Vert 
F_A\Vert^2_{L^2(M)}-\frac{{\bf i}\theta}{16\pi^2}(F_A, *F_A)_{L^2(M)}}
{\rm D}[\nabla_A].\]
For simplicity we suppose that (i) the integral above localizes to 
classical solutions; (ii) all finite energy classical
solutions are self-dual; (iii) all self-dual solutions are admissibile
in the sense of Definition \ref{feltetelek}; (iv) $M$ is simply 
connected and the surface $B_{+\infty}$ is orientable. Let $\nabla_{A_0}$ be a 
finite energy classical solution to the Yang--Mills equations. Then by 
Theorem \ref{fotetel} there exists an integer $k\in\N$ and a flat 
connection $\nabla_\Gamma\vert_W$ with $m=0$ i.e., which is a pullback 
of a flat connection $\nabla_{\Gamma_{+\infty}}$ on $B_{+\infty}$, such 
that $([\nabla_{A_0}],\Gamma )\in{\cm}(k,\Gamma )$. 
Then the previous integral cuts down to
\[Z(M,g,\tau, {\rm SU}(2))=Z_{\rm Quant}\sum\limits_{k\in\N}
\int\limits_{\cm (B_{+\infty})}\int\limits_{\widehat{\cm}(k,\Gamma )}
{\rm e}^{\frac{1}{2e^2}\Vert
F_{A_0}\Vert^2_{L^2(M)}-\frac{{\bf i}\theta}{16\pi^2}(F_{A_0},
*F_{A_0})_{L^2(M)}}{\rm D}[\nabla_{A_0}]{\rm D}
[\nabla_{\Gamma_{+\infty}}]\]
where $\cm (B_{+\infty})$ is the modulis space of flat connections over 
$B_{+\infty}$. We may try to calculate this integral within the framework of 
(topological) quantum field theory. Since elements of $\cb (E_0)$ 
extend as finite energy objects over $X$ and in particular elements of 
$\widehat{\cm} (k,\Gamma )$ remain smooth on it it is plausible to replace $M$ 
by its compactification $X$ as in (\ref{vegtelen}) and suppose that 
$Z(M,g,\tau, {\rm SU}(2))=Z(X,\tau ,{\rm SU}(2))$ where $Z(X,\tau, {\rm 
SU}(2)):\ch_{-\infty}(\emptyset )\rightarrow\ch_{+\infty}(\emptyset )$ 
is a linear map and $\ch_{\pm\infty}(\emptyset )\cong\C$ are the Hilbert 
spaces attached to the past and future boundaries of the closed space $X$ now 
considered as a cobordism between two emptysets. Assume that for a fixed 
$0<R<+\infty$ the space $X$ is cut up along 
$\partial\overline{M}_R=N\times\{ R\}$ as follows:
\[ X=(M\setminus V_R)\cup_{N\times\{ R\}} V_R\]
and let $\ch_R(N)$ denote the Hilbert space associated to 
$\partial\overline{M}_R\cong N$. By the standard axioms we expect that 
there exist vectors $v_R\in \ch_R(N)$ and $w_R\in\ch_R(N)^*$ such that 
$Z(X,\tau ,{\rm SU}(2))=(v_R, w_R)$ and the left hand side is 
independent of the particular value of $R$. Therefore taking the limit 
$R\rightarrow +\infty$ we formally obtain
\begin{equation}
Z(M,g,\tau ,{\rm SU}(2))=Z(X,\tau ,{\rm SU}(2))=(v_{+\infty}, 
w_{+\infty})
\label{particiofv}
\end{equation}
where $v_{+\infty}\in \ch_{+\infty}(N)\cong\ch (B_{+\infty})$ and 
similarly $w_{+\infty}\in\ch(B_{+\infty})^*$ since in the limit 
$R\rightarrow +\infty$ the fibration (\ref{fibralas}) cuts down to 
$B_{+\infty}\subset X$ as in (\ref{vegtelen}).

What sort of space is $\ch (B_{+\infty})$ here? We can follow the 
original ideas of Witten \cite{wit1}. Theorem \ref{fotetel} says that if $M$ 
is simply connected then all admissible 
${\rm SU}(2)$ Yang--Mills instantons approach a flat connection  
$\nabla_{\Gamma}\vert_{V_R}$ which in the limit $R\rightarrow +\infty$ 
smoothly reduces to a flat ${\rm SU}(2)$ connection 
$\nabla_{\Gamma_{+\infty}}$ on $B_{+\infty}$. Therefore by the principles of 
geometric quantization we would expect that
\[\ch(B_{+\infty})=\bigoplus_{k\in\N}H^0\left(\cm (B_{+\infty}); \co (L^k)
\right)\] 
where $L$ is the usual quantizing line bundle over $\cm (B_{+\infty})$. 
To regard $H^0$ as the space of holomorphic sections we need a complex 
structure on $\cm (B_{+\infty})$ which is inherited from one on $B_{+\infty}$. 
However the whole procedure and hence the space $\ch(B_{+\infty})$ is 
expected to be independent of any particular complex structure which leads 
to the usual conclusion that 
$H^0\left(\cm (B_{+\infty}); \co (L^k)\right)$ should be the space 
of conformal blocks of some conformal field theory (probably the 
${\rm SU}(2)$ Wess--Zumino--Witten model over $B_{+\infty}$ at level $k$). 

Therefore $\ch (B_{+\infty})$ would carry a representation of the symmetry 
algebra of some conformal field theory and in particular of the mapping 
class group of $B_{+\infty}$; these might lead to the understanding 
of the modular properties of the original partition function 
$Z(M,g,\tau ,{\rm SU}(2))$ if written in the form (\ref{particiofv}). In this 
way apparently a very geometric link between $4d$ YM theory over an ALF space 
and $2d$ CFT emerges which supports some recent investigations 
\cite{ald-gai-tac, alf-tar, bel-ber-fei-lit-tar, tac}.

\vspace{0.1in}

\noindent {\bf Acknowledgement.} This note can be considered as a 
completion of a series of papers \cite{ete,ete-hau1, ete-hau2, 
ete-hau3, ete-jar, ete-sza} on the subjectmatter. The author would like to 
thank to all of his co-authors, namely T. Hausel, M. Jardim, Sz. Szab\'o for 
the stimulating collaborations around the world during the past decade. 
Also thanks go to R.A. Mosna for finding new continuous energy irreducible 
${\rm SU}(2)$ instantons over the Schwarzschild geometry \cite{mos-tav}. 
The author also would like to thank to K. Wehrheim and the Referees of CMP in 
helping to clarify several technical points in this note. This work was 
supported by OTKA grant No. NK81203 (Hungary).

\end{document}